\documentclass[a4paper]{jpconf}
\usepackage{iwsmihead}
\usepackage{graphicx}
\usepackage{amsfonts}

\newcommand{\cbJ}{\mbox{\boldmath{$\cal J$}}}

\newcommand{\bX}{\mbox{\boldmath{$X$}}}

\newcommand{\bS}{\mbox{\boldmath{$S$}}}

\newcommand{\bv}{\mbox{\boldmath{$v$}}}

\newcommand{\bd}{\mbox{\boldmath{$d$}}}

\newcommand{\bI}{\mbox{\boldmath{$I$}}}
\newcommand{\bJ}{\mbox{\boldmath{$J$}}}

\begin{document}
\title{Cavity approach to the first eigenvalue problem in a family of symmetric random sparse matrices}

\author{Yoshiyuki Kabashima$^{\dag 1}$, Hisanao Takahashi$^{\dag 2}$ and Osamu Watanabe$^{\ddag 3}$}

\address{$^\dag$Department of Computational Intelligence and Systems Science, Tokyo Institute of Technology, Yokohama 226-8502, Japan \\
$^\ddag$Department of Mathematical and Computing Science, Tokyo Institute of Technology, Tokyo 152-8552, Japan}

\ead{$^1$kaba@dis.titech.ac.jp, $^2$takahashi@sp.dis.titech.ac.jp, $^3$watanabe@is.titech.ac.jp}

\begin{abstract}
A methodology to analyze the properties of the first (largest) eigenvalue and its eigenvector is developed for large symmetric random sparse matrices utilizing the cavity method of statistical mechanics. Under a tree approximation, which is plausible for infinitely large systems, in conjunction with the introduction of a Lagrange multiplier for constraining the length of the eigenvector, the eigenvalue problem is reduced to a bunch of optimization problems of a quadratic function of a single variable, and the coefficients of the first and the second order terms of the functions act as cavity fields that are handled in cavity analysis. We show that the first eigenvalue is determined in such a way that the distribution of the cavity fields has a finite value for the second order moment with respect to the cavity fields of the first order coefficient. The validity and utility of the developed methodology  are examined by applying it to two analytically solvable and one simple but non-trivial examples in conjunction with numerical justification. 
\end{abstract}

\section{Introduction} 
The first (largest) eigenvalue and its eigenvector (first eigenvector) play key roles in many problems in information science. In multivariate data analysis, the first eigenvector of the variance-covariance matrix represents the most significant component that underlies a set of data under the assumption that the data are generated from a multivariate Gaussian distribution, and the first eigenvalue indicates its relevance \cite{PCA}. The well-known Google PageRank$^{\rm TM}$ ranks World Wide Web pages based on the first eigenvector of a transition matrix over a huge network of pages \cite{PageRank}. The first eigenvectors of certain matrices expressing a given graph can also be utilized as a practical solution to combinatorial problems such as graph three coloring- \cite{Alon} and graph bisection problems \cite{CojaOghlan}.

The first eigenvalue problem is of significance in physics as well. The assessment of the ground state is generally formulated as a first eigenvalue problem in quantum mechanics \cite{Quantum}. It is generally difficult to accurately assess the correlations in spin systems in statistical mechanics. However, when the temperature is sufficiently high, spin correlations can be evaluated by replacing spin variables with so-called spherical spins in a class of systems \cite{ParisiPotters}, which is practically reduced to the eigenvalue analysis of the interaction matrix. In particular, that for the first eigenvalue is of great relevance because it directly leads to the assessment of the critical temperature/mode for the emergence of spontaneous magnetization \cite{FischerHertz}. 

Many properties have thus far been clarified for the first eigenvalue problem for the ensembles of large dense matrices. For $N \times N$ symmetric matrices $\bS$ whose entries are independently and identically distributed (i.i.d.) random variables of a zero mean and a variance of $N^{-1}$, the first eigenvalue asymptotically converges to $2$ and the first eigenvector has no preferential directions in the $N$-dimensional space as $N$ tends toward infinity \cite{RandomMatrix}. However, when a projection operator of a certain direction, $\bd$, is added to $\bS$ with an amplitude, $B$, the asymptotic eigenvalue is switched to $B+1/B$ for $B>1$ having a first eigenvector of a non-zero overlap with $\bd$. Similar results are also obtained for the correlation matrix, $P^{-1}\bX^{\rm T}\bX$, of $P \times N$ matrices $\bX$ whose entries i.i.d. random variables of a zero mean and variance $N^{-1}$, where ${\rm T}$ denotes the operation of the matrix transpose; the first eigenvalue is offered as $(1+\alpha^{-1/2})^2$ in the limit of $N,P \to \infty$ keeping $\alpha=P/N$ finite providing no preferential directions of the first eigenvector \cite{MarchenkoPastur}, and the analysis can also be generalized to cases in which the projection operators of several directions with arbitrary strengths are added \cite{Magnus}. However, as far as we know, relatively little knowledge has been gained on ensembles of sparse matrices in which the density of non-zero entries vanishes as the size of the matrices tends to infinity although details on eigenvalue distributions have recently been unraveled for several cases \cite{BrayRodgers,BiroliMonasson,Semerjan,NagaoTanbaka,Kuhn,Takeda}. 

Based on such an understanding, we will herein examine the first eigenvalue and the distribution of components of the first eigenvector for a family of large symmetric random sparse matrices. To achieve this purpose, we employ the cavity method of statistical mechanics \cite{Cavity}. A major advantage of this scheme is its capability for deriving equations for macroscopically characterizing objective systems without having to resort to complicated computation that is generally demanded in an alternative approach termed the replica method \cite{Replica}, on the basis of a tree approximation, which is intuitively plausible for large sparse matrices 

This paper is organized as follows. The next section introduces the model that will be explored. In section 3, we develop a scheme for examining the eigenvalue problem in a large system limit based on the cavity method. The scheme is applied to two analytically solvable and one simple but nontrivial examples in section 4. The final section is devoted to a summary. 

\section{Model definition} 
We consider ensembles of $N \times N$ real symmetric sparse matrices $\bJ=(J_{ij})$ that are characterized by a distribution, $p(k)$, of {\em degree} $k (=0,1,2,\ldots)$, which denotes the number of non-zero entries for a column/row in the limit of $N \to \infty$. For simplicity, we assume that the diagonal elements of the matrices are always constrained to zero. For aspects other than degrees, we assume that the matrices are randomly constructed. Let us denote $d_i$ as the degree of index $i(=1,2,\ldots,N$). For $k=0,1,2,\ldots$, we set $d_i=k$ for $Np(k)$ indices of $i=1,2,\ldots,N$. A practical scheme for generating a random configuration of non-zero entries characterized in the above is basically as follows \cite{StegerWormald}: 
\begin{itemize} 
\item[(S)] Make a set of indices $U$ to which each index $i$ attends $d_i$ times. Accordingly, we iterate the following (A)--(C). 
\item[(A)] Choose a pair of two different elements from $U$ randomly. 
\item[(B)] Let us denote the values of the two elements as $i$ and $j$. If $i \ne j$ and the pair of $i$ and $j$ has not been chosen up to the moment, make a link between $i$ and $j$, and remove the two elements from $U$. Otherwise, we return them back to $U$. 
\item[(C)] If $U$ becomes empty, finish the iteration. Otherwise, if there is no possibility that any more links can be made by (A) and (B), return to (S). 
\end{itemize}

Links generated by the above procedure stand for pairs of indices $i$ and $j$ to which the non-zero entries of $J_{ij}$ are assigned. There may be other algorithms for randomly generating a configuration of the non-zero entries. However, we expect that the properties, which we will investigate after this, do not depend on the details of the generation schemes. 

When the support of $p(k)$ is not bounded from the above and values of the entries are kept finite, the first eigenvalue generally diverges as $N \to \infty$ \cite{BrayRodgers,BiroliMonasson,Semerjan,NagaoTanbaka,Kuhn,Takeda}. To avoid this possibility, we assume that $p(k)=0$ for $k$, which is larger than a certain value, $k_{\rm max}$, unless infinitesimal entries are assumed. We also assume that the non-zero elements of $J_{ij}$ are determined as statistically independent samples following an identical distribution, $p_{\rm J} (J_{ij})$. In cases where $p_{\rm J}(J_{ij})$ is provided by a simple binary distribution, 

\begin{eqnarray} 
p_{\rm J}(J_{ij}|\Delta,J)= \frac{1+\Delta}{2}\delta(J_{ij}-J) + \frac{1-\Delta}{2} \delta(J_{ij}+J) 
\label{binary_dist} 
\end{eqnarray}
$(1\ge \Delta \ge 0, \ J>0)$, where $\delta(x)$ denotes Dirac's delta function, these guarantee that the first eigenvalue $\Lambda$ is bounded from the above by $k_{\rm max} J$. 

In matrix ensembles of Erd\"{o}s-R\'{e}nyi type \cite{ErdosRenyi}, which have been widely studied in network science, a non-zero value of $J_{ij}$ is assigned to each pair of indices $i>j$ with a certain probability, $c/N$, where $c>0$ is $O(1)$. This eventually leads to a Poissonian degree distribution as $p(k)=e^{-c} c^k/k !$, whose support is not bounded from the above. Therefore, the following analysis does not cover such ensembles as long as the matrix entries are $O(1)$. 

\section{Cavity approach to first eigenvalue problem} 
\subsection{Message passing in fixed graph} 
Formulating the first eigenvalue problem as a constrained quadratic optimization problem 
\begin{eqnarray} 
\mathop{\rm min}_{\bv} \left \{-\bv^{\rm T}\bJ \bv \right \} \ {\rm subject \ to} \ |\bv|^2 =N 
\label{constrained} 
\end{eqnarray}
is the basis of our analysis, where $\mbox{\rm min}_{X} \{\cdots \}$ denotes minimization with respect to $X$. The solution to this problem $\bv^*=(v_i^*)$ represents the first eigenvector that is normalized to $\sqrt{N}$ and the first eigenvalue $\Lambda$ is assessed as $\Lambda=(\bv^*)^{\rm T} \bJ \bv^*/N$. The introduction of a Lagrange multiplier $\lambda$ converts (\ref{constrained}) to a saddle point problem with respect to an objective function, 
\begin{eqnarray} 
{\cal L}(\bv,\lambda) &=& -\bv^{\rm T} \bJ \bv +\lambda (|\bv|^2-N) \cr &=& \lambda \sum_{i=1} v_i^2 - 2 \sum_{i > j} J_{ij}v_i v_j - N \lambda. 
\label{unconstrained} 
\end{eqnarray}
{Minimizing} this function with respect to $\bv$ and determining $\lambda$ so as to satisfy the stationary condition, $|\bv|^2=N$, indicate that the solution is provided by the first eigenvalue, $\lambda=\Lambda$, and the first eigenvector, $\bv=\bv^*$. 

The key idea underlying our approach is to approximate the minimization problem of (\ref{unconstrained}) with respect to $\bv$ for fixed $\lambda$ by a bunch of those for effective quadratic functions of a single variable 
\begin{eqnarray} 
{\cal L}_i(v_i|A_i,H_i)=A_i v_i^2 - 2 H_i v_i 
\label{single_quadratic_function} 
\end{eqnarray}
$(i=1,2,\ldots,N)$, where the coefficients of the second and first order terms, $A_i$ and $H_i$, are determined in a certain self-consistent manner. In the cavity method, these coefficients are evaluated by a message-passing algorithm that yields the exact solution when the connectivity of $J_{ij}$ is pictorially represented by a cycle free graph (tree). 

To explain this, let us focus on an arbitrary index, $i$, assuming that matrix connectivity is specified by a tree. We use a notation, $\partial i$, to represent the set of indices $j$ that are connected directly to $i$ with non-zero entries of $J_{ij}$. We also introduce auxiliary variables $A_{j \to i}$ and $H_{j \to i}$ to represent the second and first order coefficients of $j \in \partial i$ in the {\em $i$-cavity system}, which is defined by removing index $i$ from the original system. In physics, such variables are occasionally termed {\em cavity fields}. 

A general and distinctive feature of trees is that all indices $j \in \partial i$ are completely disjointed from one another by removing the index, $i$. Let us assume that we add $i$ into the $i$-cavity system retaining one link $J_{li}$ for an index, $l \in \partial i$, removed. This creates a problem with minimizing a function with respect to the $v_i$ and $v_j$ of $j \in \partial i \backslash l$, where $S \backslash l$ indicates the removal of $l$ from a set, $S$, as 
\begin{eqnarray} 
{\cal L}_{i, \partial i \backslash l}(v_i,\{v_{j \in \partial i \backslash l}\}) =\lambda v_i^2 -2 v_i \sum_{j \in \partial i \backslash l}J_{ij} v_j + \sum_{j \in \partial i \backslash l} \left (A_{j \to i}v_j^2 - 2 H_{j \to i} v_j \right ). 
\label{cavity_update_pre} 
\end{eqnarray}
Since $i$ is connected to $l$ only with $J_{li}$ in the original tree, minimizing this function with respect to the $v_j$ of $\forall{j} \in \partial i \backslash l$ fixing $v_i$ yields cavity fields concerning $i$ in the $l$-cavity system as 
\begin{eqnarray} 
&&A_{i \to l} = \lambda- \sum_{j \in \partial i \backslash l} \frac{J_{ij}^2}{A_{j \to i}}, 
\label{cavity2} \\ && H_{i \to l} = \sum_{j \in \partial i \backslash l} \frac{J_{ij} H_{j \to i}}{A_{j \to i}}. 
\label{cavity1} 
\end{eqnarray}
Given an initial condition, these equations are capable of assessing sets of the cavity fields, $A_{j \to i}$ and $H_{j \to i}$, which are defined for all directed pairs of the connected indices, moving entirely over the graph. $J_{ij}^2/A_{j \to i} $ and $J_{ij} H_{j \to i}/A_{j \to i} $ in (\ref{cavity2}) and (\ref{cavity1}) represent the influences of $j$ to $i$ conveyed through connection $J_{ij}$ and are sometimes referred to as {\em cavity biases}. After the cavity fields are determined, $A_i$ and $H_i$ in (\ref{single_quadratic_function}) are assessed by summing up the cavity biases from $\forall{j} \in \partial i$ as $A_i = \lambda - \sum_{j \in \partial i} {J_{ij}^2}/{A_{j \to i}}$ and $H_i = \sum_{j \in \partial i} J_{ij} H_{j \to i}/A_{j \to i}$. This yields $v_i^*=H_i/A_i$. Finally, adjusting $\lambda$ so that the normalization constraint, $|\bv^*|^2 =\sum_{i=1}^N (v_i^*)^2=\sum_{i=1}^N (H_i/A_i)^2=N$, is satisfied provides $\Lambda$. 

The above procedure can be regarded as a variant of the belief propagation developed in research on probabilistic inference \cite{Pearl} and is guaranteed to offer the exact result when the graph is free from cycles. When the graph contains cycles, the algorithm in (\ref{cavity2}) and (\ref{cavity1}) can still be performed; unfortunately, the result obtained is just an approximation. However, for our current purpose, we can generally utilize a much simpler and exact method that just repeats matrix multiplication as $\bv^{t+1}= (a \bI_N+\bJ) \bv^t$, 
which provides $\Lambda=\lim_{t \to \infty} |\bv^{t+1}|/|\bv^t|-a $ and $\bv^*=\sqrt{N} \lim_{t \to \infty} \bv^t/|\bv^t|$, where $a$ is a sufficiently large positive number and $\bI_N$ denotes the $N \times N$ identity matrix.  
Therefore, our approach may seem less competitive in practice. However, the cavity approach is still useful for examining the properties of the eigenvalue problem introducing a macroscopic description as is explained below. 

\subsection{Description by distribution of cavity fields} 
Let us examine the behavior of the algorithm in (\ref{cavity2}) and (\ref{cavity1}). To do this, we need to pay attention to the property of randomly constructed sparse matrices, which implies that the lengths of cycles in the connectivity graph of a matrix typically grow as $O(\ln N)$ as $N \to \infty$ \cite{RandomGraph}. This naturally motivates us to ignore the effects of self-interactions in the updates of (\ref{cavity2}) and (\ref{cavity1}). We also introduce a macroscopic characterization of the cavity fields utilizing the distribution, $q(A,H)=(\sum_{i=1}^N d_i)^{-1} \sum_{i=1}^N \sum_{j \in \partial i} \delta(A-A_{j \to i}) \delta(H-H_{j \to i})$. When a link in the connectivity graph is chosen randomly, the probability that the index of one terminal has a degree, $k$, is provided as 
\begin{eqnarray} 
r(k)=\frac{k p(k)}{\sum_{k=0}^{k_{\rm max}} k p(k)}. 
\label{posterior} 
\end{eqnarray}
This and dealing with (\ref{cavity2}) and (\ref{cavity1}) as the elemental process for updating the cavity field distribution lead to an equation that determines $q(A,H)$ in a self-consistent manner as 
\begin{eqnarray} 
q(A,H) \! =\! \sum_{k=1}^{k_{\rm max}} \! r(k) \! \int \! \prod_{j=1}^{k-1} \! dA_j dH_j q(A_j,H_j) \left \langle \delta\left (A-\lambda+\sum_{j=1}^{k-1}\frac{{\cal J}_j^2}{A_j} \right ) \delta\left (H-\sum_{j=1}^{k-1} \frac{{\cal J}_j H_j}{A_j} \right ) \right \rangle_{\cbJ} \! , 
\label{cavity_dist_update} 
\end{eqnarray}
where $\left \langle \cdots \right \rangle_{\cbJ}$ represents the operation of taking averages with respect to all relevant ${\cal J}_j$'s following $p_{\rm J}({\cal J}_j)$. After $q(A,H)$ is determined with this equation, the distribution of the auxiliary variables in the original system, $Q(A,H)=N^{-1} \sum_{i=1}^N \delta(A-A_i)\delta(H-H_i)$, is assessed as 
\begin{eqnarray} 
Q(A,H) \! = \! \sum_{k=0}^{k_{\rm max}} \! p(k) \! \int \! \prod_{j=1}^k \! dA_j dH_j q(A_j,H_j) \left \langle \delta\left (A-\lambda+\sum_{j=1}^{k}\frac{{\cal J}_j^2}{A_j} \right ) \delta\left (H-\sum_{j=1}^{k} \frac{{\cal J}_j H_j}{A_j} \right ) \right \rangle_{\cbJ} \!. 
\label{fulldist} 
\end{eqnarray}
This make it possible to evaluate the quadratic norm per element of $\bv^*$ as $T=N^{-1}|\bv^*|^2 =\int dA dH Q(A,H) (H/A)^2$. Moreover, the distribution of elements of the first eigenvector, $\bv^*$, can be assessed as $P(v)=N^{-1}\sum_{i=1}^N \delta(v-v^*)= \int dA dH Q(A,H) \delta(v-H/A)$. 

The argument provided in the preceding subsection implies that $\lambda$ should be adjusted so that $T$ 
accords 
with unity. However, here we adopt an alternative approach to determining $\Lambda$. Let us assume a situation where (\ref{cavity_dist_update}) is solved by a method of successive iteration, which leads to an exact result when the graph is free from cycles. Equation (\ref{cavity_dist_update}) guarantees that for large $\lambda$, $A$ is distributed over large values in $q(A,H)$. This means that for sufficiently large $\lambda$, the absolute values of $H$ are mostly reduced by each iteration and, therefore, the marginal distribution, $q(H)=\int dA q(A,H)$, converges to $\delta(H)$ implying $T \to 0$ as the number of iterations tends to infinity. However, as $\lambda$ is decreased from sufficiently large values, $A$ can take smaller values in the distribution and larger values of $|H|$ could appear with considerable frequency in each iteration. This indicates that the frequencies with which $|H|$ is reduced and enlarged are balanced at a certain value of $\lambda$, for which non-trivial distribution $q(H) \ne \delta(H)$ becomes invariant under the cavity iteration of (\ref{cavity_dist_update}) and $T$ is kept as a finite constant. If $\lambda$ is lowered further, the frequency of enlargements will overcome that of reductions, which will eventually make $T$ tend to infinity. These imply that one can characterize the first eigenvalue, $\Lambda$, as the value of $\lambda$ for which non-trivial distribution $q(H)\ne \delta(H)$ emerges by solving (\ref{cavity_dist_update}). The resulting cavity field distribution, $q(A,H)$, offers the distribution of the first eigenvector elements, $P(v)$, of any finite value of $T$ by applying appropriate rescaling with respect to $H$. This is also reasonable because eigenvalues are generally irrelevant to the values of the normalization constraint with respect to eigenvectors. 

The validity and utility of this characterization of the eigenvalue problem are examined by applying it to three examples in the next section. 

\section{Application to three examples} 
\subsection{Single degree model} 
In the first example, we consider cases where all indices have a certain identical degree, $K$, which is characterized by $p(k)=\delta_{k,K}$, where $\delta_{i,j}$ denotes Kroeneker's delta, and the non-zero entries are i.i.d. following the binary distribution of (\ref{binary_dist}). In such cases, the marginal distribution, $q(A)=\int dH q(A,H)$, is provided in the form of $q(A)=\delta(A-A^*)$ since ${\cal J}^2=J^2$ always holds for $\forall{\cal J}$ that are sampled from (\ref{binary_dist}). Equation (\ref{cavity_dist_update}) indicates that $A^*$ satisfies $A^*=\lambda-(K-1)J^2/A^*$. This yields 
\begin{eqnarray} 
A^*=\frac{\lambda+\sqrt{\lambda^2-4(K-1)J^2}}{2} \equiv A^*(\lambda;(K-1)J^2), 
\label{A_first_example} 
\end{eqnarray}
which along with (\ref{cavity_dist_update}) offer the equation for determining the marginal distribution, $q(H)$, in a self-consistent manner as 

\begin{eqnarray} 
q(H)=\int \prod_{j=1}^{K-1} dH_j q(H_j) \left \langle \delta \left (H-\sum_{j=1}^{K-1} \frac{{\cal J}_j H_j}{A^*(\lambda;(K-1)J^2)} \right ) \right \rangle_{\cbJ}. 
\label{marginal_q_dist} 
\end{eqnarray}

\begin{figure}[t] 
\setlength\unitlength{1mm} 
\begin{picture}(140,60)(0,0) \put(10,-35){\includegraphics[width=14cm]{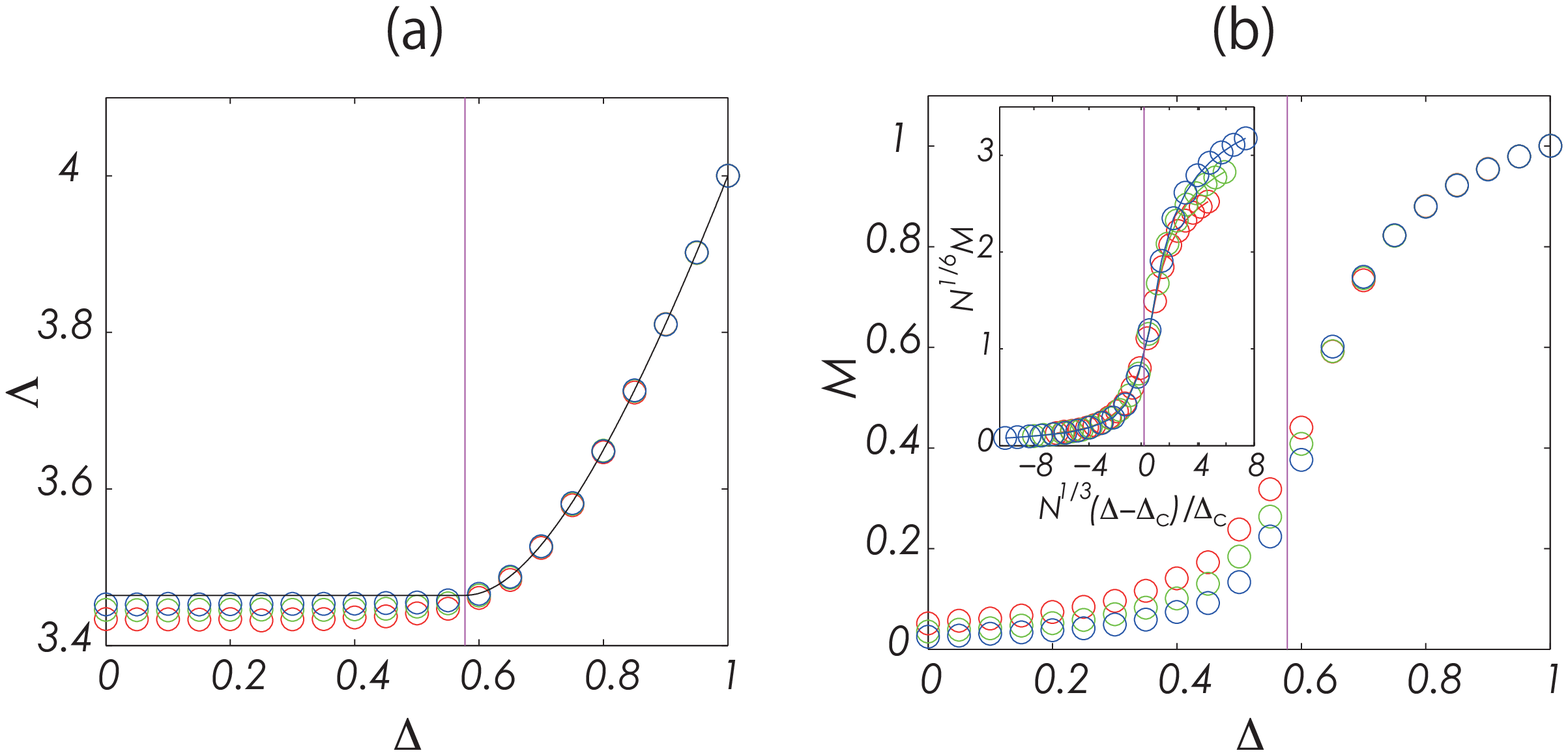}} 
\end{picture} 
\caption{(a): First eigenvalue $\Lambda$ versus $\Delta$ for single degree model of $K=4$ and $J=1$. The full curve represents the theoretical prediction of (\ref{first_eigenvalue}) while the markers denote the averages of $2000$ experimental results for $N=256$ (red), $512$ (green) and $1024$ (blue) systems. (b): $M=|N^{-1}\sum_{i=1}^N v_i^*|$ versus $\Delta$ for the same experiments as $(a)$. Red, green, and blue markers correspond to $N=256$, $512$, and $1024$ as well. The inset has scaling plots obtained by conversion of $M \to N^{1/6} M$ and $\Delta \to N^{1/3} (\Delta-\Delta_{\rm c})/\Delta_{\rm c}$. In both (a) and (b), the vertical straight lines (magenta) stand for $\Delta_{\rm c}=1/\sqrt{3}=0.577\ldots$.} 
\label{fig1} 
\end{figure}

Tuning $\lambda$ so that a non-trivial distribution, $q(H) \ne \delta(H)$, satisfies this equation yields the first eigenvalue, $\Lambda$. Paying attention to the first and the second order moments, $m_1=\int dH q(H) H$ and $m_2=\int dH q(H) H^2$, is sufficient for this purpose. This provides a generalized eigenvalue problem described by a couple of equations 
\begin{eqnarray} 
m_1&=&\frac{\Delta (K-1)J}{A^*(\lambda;(K-1)J^2)} m_1, 
\label{moment1} \\ m_2&=&\frac{(K-1)J^2}{(A^*(\lambda;(K-1)J^2))^2} (m_2-m_1^2) + \left (\frac{\Delta (K-1)J}{A^*(\lambda;(K-1)J^2)} m_1 \right )^2, 
\label{moment2} 
\end{eqnarray}
which requires the existence of a solution to $(m_1,m_2) \ne (0,0)$. There are two possibilities that will satisfy this requirement. The first is characterized by $|m_1| > 0$, which in conjunction with (\ref{moment1}) offers 
\begin{eqnarray} 
1=\frac{\Delta (K-1)J}{A^*(\lambda;(K-1)J^2)}. 
\label{first_moment} 
\end{eqnarray}
The second is the case of $m_1=0$ and $m_2 > 0$, which along with (\ref{moment2}) provide 
\begin{eqnarray} 
1=\frac{(K-1)J^2}{(A^*(\lambda;(K-1)J^2))^2}. 
\label{second_moment} 
\end{eqnarray}
Solving (\ref{first_moment}) and (\ref{second_moment}) with respect to $\lambda$, we finally obtain an expression of the first eigenvalue, 
\begin{eqnarray} 
\Lambda=\left \{
\begin{array}{ll} \left ((K-1)\Delta+1/\Delta \right )J, & \Delta>\Delta_{\rm c}, \cr 2 \sqrt{K-1} J, & \Delta < \Delta_{\rm c}, 
\end{array} \right . 
\label{first_eigenvalue} 
\end{eqnarray}
where $\Delta_{\rm c}=1/\sqrt{K-1}$. $\Delta > \Delta_{\rm c}$ and $\Delta < \Delta_{\rm c}$ correspond to the cases of $|m_1| >0$ and $m_1 =0$. 

\begin{figure}[t] 
\setlength\unitlength{1mm} 
\begin{picture}(140,60)(0,0) \put(5,-20){\includegraphics[width=14cm]{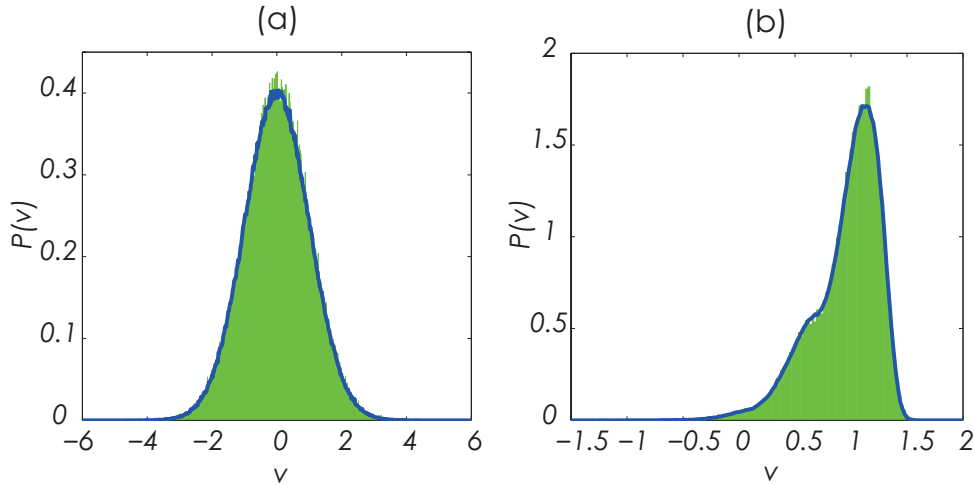}} 
\end{picture} 
\caption{$P(v)$ for (a): $\Delta=0.5$ and (b): $\Delta=0.9$ of single degree model of $K=4$. For both plots, the bars represent histograms obtained from $100$ experiments of $N=1024$ while the curves stand for theoretical predictions obtained from (\ref{marginal_q_dist}) and (\ref{fulldist}) by means of a Monte Carlo method of $10^6$ populations. In handling the experimental data of (b), we employed conversion $\bv^* \to -\bv^*$ for data of $\sum_{i=1}^N v_i^* < 0$ in order to break the mirror symmetry between $\bv^*$ and $-\bv^*$, which is intrinsic in the eigenvalue problem.} 
\label{fig2} 
\end{figure}

To justify the prediction of (\ref{first_eigenvalue}), we carried out numerical experiments for the systems of $K=4$ and $J=1$, whose results are in figures \ref{fig1} (a) and (b). Figure \ref{fig1} (a) has the plots of $\Lambda$ versus $\Delta$, which are in excellent agreement with (\ref{first_eigenvalue}) for relatively large $\Delta$. However, there are non-negligible finite-size corrections for relatively small $\Delta$, which makes it difficult to accurately detect the critical value, $\Delta_c$. To overcome this difficulty, we turn to the profiles of $M=\left |N^{-1} \sum_{i=1}^N v_i^* \right |$, which are shown in 
figure \ref{fig1} (b). The raw experimental data of $M$ vary smoothly with $\Delta$ and do not exhibit any singularity at $\Delta=\Delta_{\rm c}$ either. Nevertheless, the data rescaled as $M \to N^{1/6} M$ and $\Delta \to N^{1/3} (\Delta-\Delta_{\rm c})/\Delta_{\rm c}$ support a scaling hypothesis, $M=N^{-1/6} F(N^{1/3}(\Delta-\Delta_{\rm c})/\Delta_{\rm c})$ (inset), which validates the transition at $\Delta_{\rm c}$ assuming that $F(x)$ behaves as $O( x^{1/2})$ for $x \to +\infty$ and vanishes for $x \to -\infty$. A similar scaling hypothesis has been assumed in the examination of principal component analysis before \cite{Magnus}.

Besides the eigenvalue, one can also assess the distributions, $P(v)$, of elements of the first eigenvector, $\bv^*$, by numerically solving $q(H)$ with Monte Carlo (population dynamics) methods. The prediction with our approach is also in excellent agreement with the experimental results for both cases of $\Delta < \Delta_{\rm c}$ (figure \ref{fig2} (a)) and $\Delta > \Delta_{\rm c}$ (figure \ref{fig2} (b)) when $|\Delta -\Delta_{\rm c}|/|\Delta_{\rm c}|$ is sufficiently large. 

\subsection{Limit of large degrees and infinitesimal entries} 
The second example is offered by assuming that the average and the variance of $p(k)$ are provided as $\overline{k}$ and $O(\overline{k})$, respectively, and those of $p_{\rm J}(J_{ij})$ are given as $\mu J/\overline{k}$ and $J^2/\overline{k}$ for sufficiently large $\overline{k}$ while $\mu > 0$ and $J>0$ are certain finite constants. The assumption on the degree distribution, $p(k)$, covers the cases described by truncated Poissonian distributions of large mean values in which indices whose degrees are sufficiently larger than the means are eliminated from the basic graphs of Erd\"{o}s-R\'{e}nyi type \cite{CojaOghlan}. Being combined with the assumption on $p_{\rm J}(J_{ji})$, this, in conjunction with the law of large numbers, guarantees that the marginal distribution of $A$ is provided as $q(A)=\delta(A-A^*(\lambda;J^2))$, where the functional form of $A^*(\lambda;J^2)$ is identical to that of (\ref{A_first_example}). Inserting this into (\ref{marginal_q_dist}) and the resulting equations with respect to the first and the second order moments 
offer conditions of $1=\mu J/A^*(\lambda;J^2)$ and $1=J^2/(A^*(\lambda;J^2))^2$, respectively, which are counterparts of (\ref{first_moment}) and (\ref{second_moment}). This provides the first eigenvalue as 
\begin{eqnarray} 
\Lambda=\left \{
\begin{array}{ll} (\mu + 1/\mu) J, & \mu > 1, \cr 2 J, & \mu < 1. 
\end{array} \right . 
\label{dense_eigenvalue} 
\end{eqnarray}
$\mu>1$ and $\mu < 1$ correspond to situations of $|m_1| > 0$ and $m_1=0$. 

In the current case, the central limit theorem guarantees that $q(H)$ converges to a distribution of the Gaussian type. Consequently, this also makes it possible to analytically express the distribution, $P(v)$, of the elements of the first eigenvector in a Gaussian form. Under the constraint of $T=1$, this yields 
\begin{eqnarray} 
P(v)=\left \{
\begin{array}{ll} (2 \pi (1-(M(\mu))^2)^{-1/2} \exp \left (-(v-M(\mu))^2/(2(1-(M(\mu))^2)) \right ), & \mu > 1, \cr (2 \pi)^{-1/2} \exp \left (-v^2/2 \right ), & \mu < 1, 
\end{array} \right . 
\label{dense_v_dist} 
\end{eqnarray}
where $M(\mu)=\pm \sqrt{1-1/\mu^2}$. 

\subsection{Mixture of multiple degrees} 
The strategy we employed for analyzing the preceding two examples is summarized as follows. We first evaluate the marginal distribution, $q(A)=\int dH q(A,H)$, of the second order cavity fields in the form of a single delta distribution, $q(A)=\delta(A-A^*)$. This enables us to deal with $A$ as a constant, $A^*$, in assessing the marginal distribution, $q(H)=\int dA q(A,H)$, of the first order cavity fields. Consequently, the original eigenvalue problem is reformulated as a generalized eigenvalue problem with respect to $q(H)$. This imposes certain conditions on the second order cavity field, $A^*$, which is a function of $\lambda$. The first eigenvalue is evaluated by solving the conditions with respect to $\lambda$. 

\begin{figure}[t] \setlength\unitlength{1mm} 
\begin{picture}(140,60)(0,0) \put(10,-35){\includegraphics[width=14cm]{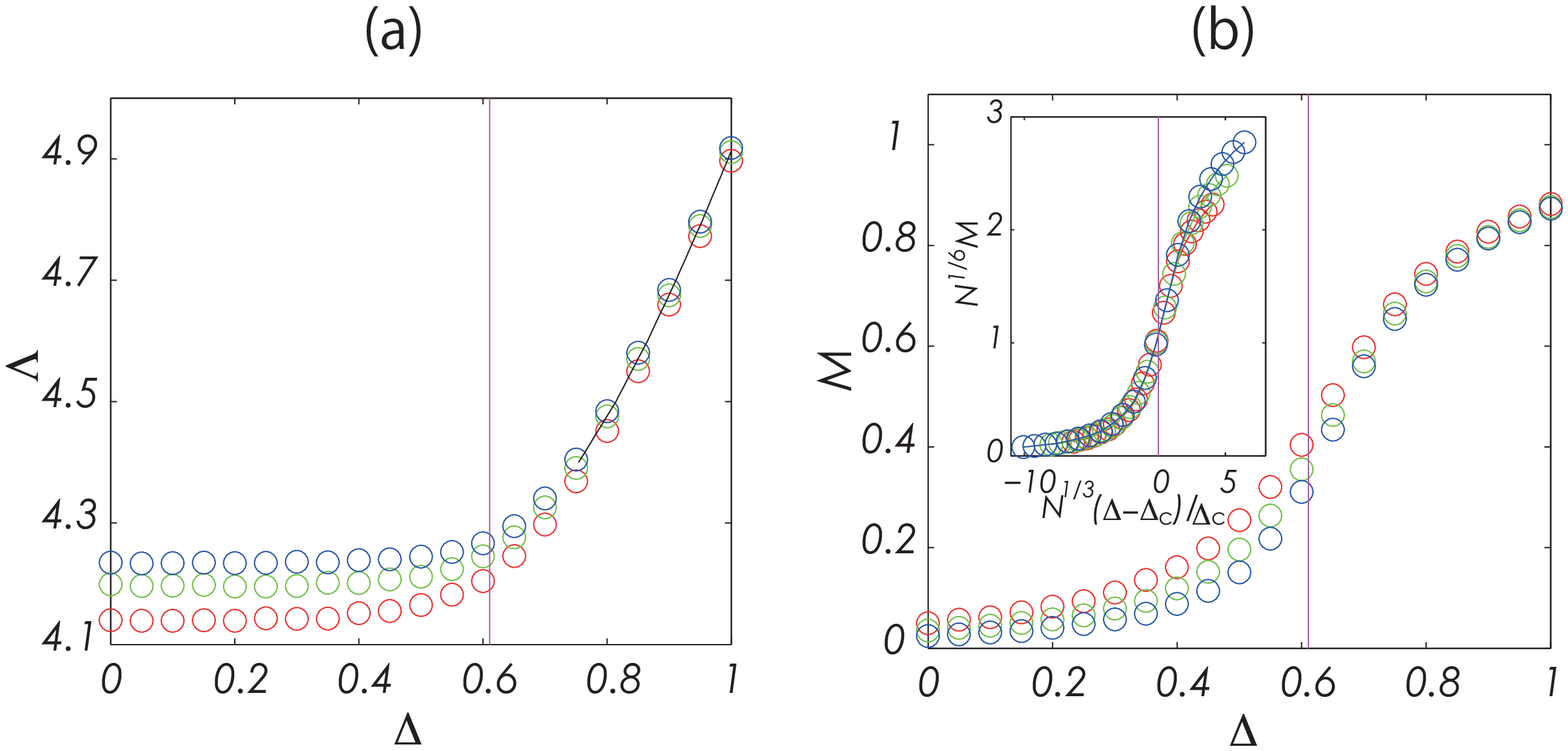}} 
\end{picture} 
\caption{(a): First eigenvalue $\Lambda$ versus $\Delta$ for mixture model of $p(k)=0.9 \delta_{k,4}+0.1 \delta_{k,8}$ and $J=1$. The full curve represents the theoretical prediction assessed from (\ref{Delta_Lambda}) by means of a Monte Carlo method of $10^6$ populations. Large statistical fluctuations due to a small portion of samples of $A \sim 0$ prevented us from accurately assessing the relation below $\Lambda \sim 4.4$, for which data have not been shown. The markers stand for data evaluated similarly for figure \ref{fig1}. (b): $M=|N^{-1}\sum_{i=1}^N v_i^*|$ versus $\Delta$ for the same experiments as (a). The scaling plots of $M \to N^{1/6} M$ and $\Delta \to N^{1/3} (\Delta-\Delta_{\rm c})/\Delta_{\rm c}$ (inset) indicate $\Delta_{\rm c} \simeq 0.611$, which is represented as vertical straight lines (magenta) in (a) and (b).} 
\label{fig3} 
\end{figure}

Unfortunately, such clear decoupling of the first and the second order cavity fields cannot generally be exploited in handling the eigenvalue problem precisely. Even when $p(k)$ takes finite values for multiple $k$, the marginal distribution, $q(A)$, constitutes a closed equation with no necessity for considering the first order cavity field, $H$, as 

\begin{eqnarray} 
q(A)=\sum_{k=1}^{k_{\rm max}}r(k) \int \prod_{j=1}^{k-1} dA_j q(A_j) \left \langle \delta \left (A-\lambda+ \sum_{j=1}^{k-1} \frac{{\cal J}_j^2}{A_j} \right ) \right \rangle_{\cbJ}, 
\label{marginal_A} 
\end{eqnarray}
which, however, generally makes $q(A)$ a continuous distribution. This prevents us from handling $A$ as a constant in assessing the marginal distribution, $q(H)$, and we eventually have to deal with the joint distribution, $q(A,H)$, directly. This requires much more effort in analysis than that in handling the marginal distributions, $q(A)$ and $q(H)$, separately. 

To avoid this, let us move forward with our analysis, approximately, ignoring the correlations between $A$ and $H$. Applying this approximate treatment to the assessment of the first and the second order moments of $H$ by using (\ref{cavity_dist_update}) provides two conditions 
\begin{eqnarray} 
1=\left (\int d{\cal J} p_{\rm J} ({\cal J}) {\cal J} \right )
\left (\sum_{k=1}^{k_{\rm max}} r(k) (k-1) \right ) \left (\int dA q(A)A^{-1} \right ), 
\label{dist1} \\ 
1=\left (\int d{\cal J} p_{\rm J} ({\cal J}) {\cal J}^2 \right ) \left (\sum_{k=1}^{k_{\rm max}} r(k) (k-1) \right ) \left (\int dA q(A) A^{-2} \right ), 
\label{dist2} 
\end{eqnarray}
where $q(A)$ is the solution to (\ref{marginal_A}) for given $\lambda$. Equations (\ref{dist1}) and (\ref{dist2}) can be regarded as corresponding to generalizations of (\ref{first_moment}) and (\ref{second_moment}). 

We examined the utility of (\ref{dist1}) and (\ref{dist2}) by applying them to cases of the bimodal degree distributions of $p(k)=(1-f) \delta_{k,4} + f \delta_{k,8}$ $(0<f<1)$ and $p_{\rm J}(J_{ij})$ of (\ref{binary_dist}) with $J=1$. 
For (\ref{binary_dist}), equation (\ref{dist1}) indicates that 
\begin{eqnarray} 
\Delta=J^{-1} \left ( \left (\sum_{k=1}^{k_{\rm max}} r(k) (k-1) \right ) \left ( \int dA q(A) A^{-1} \right ) \right )^{-1}, 
\label{Delta_Lambda} 
\end{eqnarray}
provides the relation between first eigenvalue $\Lambda$ and $\Delta$ by dealing with $\lambda=\Lambda$ as a control parameter, as long as the right hand side is less than unity and greater than a certain critical value, $\Delta_{\rm c}$. Figure \ref{fig3} (a) compares the theoretical prediction of (\ref{Delta_Lambda}) assessed with the Monte Carlo method and the results of numerical experiments for $f=0.1$. Despite the fact that the treatment is not necessarily exact, the theoretical prediction agrees with the experimental results with excellent accuracy for relatively large $\Delta$. 

\begin{figure}[t] \setlength\unitlength{1mm} 
\begin{picture}(140,60)(0,0) \put(5,-35){\includegraphics[width=14cm]{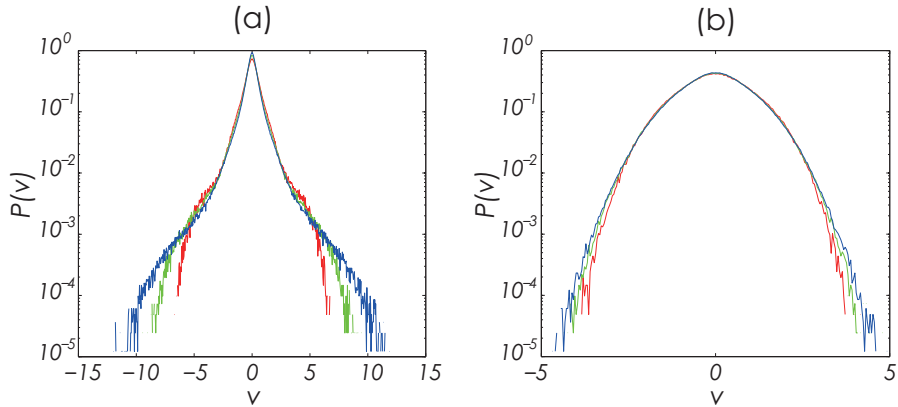}} 
\end{picture} 
\caption{$P(v)$ of mixture model $p(k)=(1-f)\delta_{k,4}+f \delta_{k,8}$ for (a): $f=0.1$ and (b): $f=0.9$. In both plots, the red, green and blue curves were obtained from $2000$ experiments setting $\Delta=0$ for $N=256$, $512$ and $1024$, respectively.} 
\label{fig4} 
\end{figure}

Rescaling $M=|N^{-1} \sum_{i=1}^N v_i^*|$ and $(\Delta-\Delta_{\rm c})/\Delta_{\rm c}$ with $N^{1/6}$ and $N^{1/3}$, respectively, as has been assumed in the single degree model for the experimental data, indicates a critical value of $\Delta_{\rm c} \simeq 0.611$ (figure \ref{fig3} (b)). This should be characterized by satisfying (\ref{dist2}) in the current approach. However, it is difficult to accurately determine $\Delta_{\rm c}$ with the Monte Carlo method in practice because the right hand side of (\ref{dist2}) greatly depends on the diverging contribution from a small portion of $A \sim 0$ at criticality, which is considerably sensitive to statistical fluctuations in the sampling. Another source of the difficulty is that (\ref{marginal_A}) can make the support of $q(A)$ spread to a region of {\em negative} $A$ for relatively small $\lambda$. The negative values of $A$ imply that for certain indices $i$, (\ref{single_quadratic_function}) is minimized at $v_i = \pm \infty$, which may correspond to large elements of $\bv^*$ that diverge as $N$ tends to infinity. Such elements have been argued to be {\em defects} in several earlier studies before \cite{BiroliMonasson,Semerjan,Kuhn}. The experimental data indicate that the divergent behavior of the small fraction of the first eigenvector elements, which may be related to the defects, is more significant as $f$ is relatively {\em smaller} (figures \ref{fig4} (a) and (b)). Refining the current analysis so that it can accurately handle such divergent elements is currently under way. 

\section{Summary} 
In summary, we developed a methodology for analyzing the first eigenvalue problem of large symmetric random sparse matrices utilizing the cavity method of statistical mechanics. The scheme we developed makes it possible to assess the first eigenvalue and the distribution of elements of the first eigenvector on the basis of a functional equation concerning a joint distribution with respect to two kinds of cavity fields. Its validity was tested and confirmed by using two analytically solvable examples. Unfortunately, directly employing the approach we developed is technically difficult in general cases where the matrices are characterized by multiple degrees. However, we demonstrated that an approximate treatment that ignores certain correlations between the two kinds of cavity fields offers excellent capabilities for assessing the first eigenvalue when each entry of the matrices is i.i.d. following a distribution that has a sufficiently large positive mean value. However, singular behavior can be observed in the cavity-field distribution when the average of the matrix entries is not a sufficiently large positive value. Comparison with experimental results suggested that this may be related to the divergent elements of the first eigenvector, which have been argued in earlier studies. 

In some problems, not only the first but also the second eigenvalue is of utility \cite{Expander}. In addition to making refinements to the current study, analyzing the $n$-th eigenvalue $(n=2,3,\ldots)$ by utilizing the cavity method may be an interesting project for future work. 

\ack 
This work was partially supported by Grants-in-Aid for Scientific Research on the Priority Area ``Deepening and Expansion of Statistical Mechanical Informatics'' from the Ministry of Education, Culture, Sports, Science and Technology, Japan and the JSPS Global COE program, ``Computationism as a Foundation for the Sciences'' (YK and OW).


\section*{References} 

\begin{thebibliography}{99}

\bibitem{PCA}
Mardia K V, Kent J T and Bibby J M 1979 {\em Multivariate Analysis} (London: Academic Press)

\bibitem{PageRank} 
Langville A and Meyer C 2006 {\em Google's PageRank and Beyond: The Science of Search Engine Rankings} (Princeton: Princeton University Press)

\bibitem{Alon} 
Alon N and Kahale N 1997 {\em SIAM J.\ Comput.\/} {\bf 26} 1733

\bibitem{CojaOghlan} 
Coja-Oghlan A 2006 {\em Random Struct. Algorithms} {\bf 29} 351

\bibitem{Quantum} 
Born M, Heisenberg W and Jordan P 1925 {\em Zeit. Phys.} {\bf 35} 557

\bibitem{ParisiPotters} 
Parisi G and Potters M 1995 \JPA {\bf 28} 526

\bibitem{FischerHertz} 
Fischer K H and Hertz J A 1991 {\em Spin Glasses} (Cambridge: Cambridge University Press)

\bibitem{RandomMatrix} 
Mehta M L 2004 {\em Random matrices, Third edition. Pure and Applied Mathematics 142} (Amsterdam: Elsevier/Academic Press)

\bibitem{MarchenkoPastur} 
Marchenko V A and Pastur L A 1967 {\em Mat. Sb.} {\bf 72} 507 

\bibitem{Magnus} 
Hoyle D C and Rattray M 2004 \PR E {\bf 69} 026124 

\bibitem{BrayRodgers} 
Bray A J and Rodgers G J 1988 \PR B {\bf 38} 11461

\bibitem{BiroliMonasson} 
Biroli G and Monasson R 1999 \JPA {\bf 32} L255

\bibitem{Semerjan} 
Semerjian G and Cugliandolo L F 2002 \JPA {\bf 35} 4837

\bibitem{NagaoTanbaka} 
Nagao T and Tanaka T 2007 \JPA {\bf 40} 4973

\bibitem{Kuhn} 
K\"{u}hn R 2008 \JPA {\bf 41} 295002

\bibitem{Takeda} 
Rogers T, R\'{e}rez-Castillo I, K\"{u}hn R and Takeda K 2008 \PR E {\bf 78} 031116

\bibitem{Cavity} 
M\'{e}zard M, Virasolo M A and Parisi G 1986 {\em Spin Glass Theory and Beyond} (Singapore: World Scientific)

\bibitem{Replica}
Dotzenko V S 2001 {\em Introduction to the Replica Theory of Disordered Statistical Systems} (Cambridge: Cambridge University Press) 

\bibitem{StegerWormald} 
Steger A and Wormald N C 1999 {\em Combinatorics, Probability and Computing} {\bf 8} 377

\bibitem{ErdosRenyi} 
Erd\"{o}s P and R\'{e}nyi A 1959 {\em Publicationes Mathematicae} {\bf 6} 290

\bibitem{Pearl} 
Pearl J 1988 {\em Probabilistic Reasoning in Intelligent Systems: Networks of Plausible Inference} (San Francisco: Morgan Kaufmann) 

\bibitem{RandomGraph} 
Albert R and Barab\'{a}si A L 2002 {\em Rev. Mod. Phys.} {\bf 74} 47

\bibitem{Expander} 
Hoory S, Linial N and Widgerson A 2006 {\em Bulletin (New series) of the American Mathematical Society} {\bf 43} 439 
\end{thebibliography}
\end{document}